\documentclass[10pt]{article}

\usepackage{a4wide}
\usepackage{amssymb}
\usepackage{amsfonts}
\usepackage{amsmath}
\input xy
\xyoption{arrow} \xyoption{matrix}

\date{}

\newtheorem{proposition}{Proposition}[section]
\newtheorem{theorem}[proposition]{Theorem}
\newtheorem{lemma}[proposition]{Lemma}

\newtheorem{corollary}[proposition]{Corollary}

\def\Hom{{\rm Hom}}
\def\der{\partial }

\def\nFM0{{\nu }_{F,M_0}}
\def\nFN0{{\nu }_{F,N_0}}
\def\nGN0{{\nu }_{G,N_0}}

\def\N0{ {\bf N}_0 }

\def\t{\otimes}
\def\g{\gamma}
\def\v{\varphi}
\def\ra{\rightarrow}

\def\lra{\leftrightarrow}
\def\Xpm{X^{\pm }}

\def\s{\sigma}
\def\Z{\mathbb{Z}}

\def\l1{{\lambda}_1}

\def\a{\alpha}
\def\a0{ {\alpha }_0}
\def\a1{ {\alpha }_1}

\def\l{\lambda}
\def\o{\omega}

\def\nFGM0{{\nu }_{F,G,M_0}}

\def\nFN0{{\nu}_{F,N_0}}

\def\sm{{\sigma}^m}

\def\sm1{{\sigma}^{-1}}

\def\smtp1{{\sigma}^{-t+1}}

\def\o{\omega }
\def\S1{S^{-1}}

\def\Xpm1{X^{\pm 1}_1}

\def\sPM1{{\sigma }^{\pm 1}}
\def\sMP1{{\sigma }^{\mp 1 }}

\def\d{\delta}

\def\di{{\rm d.ind}}

\def\L{\Lambda}

\def\Ytm1{Y^{t-1}}
\def\Yim1{Y^{i-1}}

\def\CK{{\cal K}}

\def\Aut{{\rm Aut}}

\def\dim{{\rm dim }}

\def\ker{ {\rm ker } }

\def\SL2Z{ {\rm SL}_2({\bf Z}) }

\def\th{ \theta }

\def\Gp1{ G^{1 , 1 } }
\def\P11{ P^{-1 , 1 } }
\def\Pp1{ P^{1 , 1 } }

\def\th{\theta}

\def\CV{{\cal V}}

\def\nCLsr{{}^\nu\kern-2pt {\cal L}^{\sigma , \rho  }}
\def\nP{{}^\nu \kern-2pt P}
\def\nL{{}^\nu\kern-2pt L}
\def\nLL{{}^\nu\kern-2pt \Lambda}
\def\nPsr{{}^\nu\kern-2pt P^{\sigma , \rho  }}
\def\nLsr{{}^\nu\kern-2pt L^{\sigma , \rho  }}
\def\nuCL{{}^\nu\kern-2pt  {\cal L}}
\def\nCLsr{{}^\nu\kern-2pt {\cal L}^{\sigma , \rho  }}
\def\nCL1m{{}^\nu\kern-2pt {\cal L}^{-1 , 1  }}
\def\x1nu{x^\frac{1}{\nu}}
\def\xm1nu{x^{-\frac{1}{\nu}}}

\def\ra{\rightarrow }

\def\CI{{\cal I}}

\def\coker{{\rm coker}}

\def\CC{ {\cal C}}

\def\nAM0{{\nu }_{{\cal A},M_0}}
\def\nAN0{{\nu }_{{\cal A},N_0}}

\def\End{ {\rm End }}

\def\det{ {\rm det }}

\def\ga{\mathfrak{a}}

\def\gp{\mathfrak{p}}

\def\GL{{\rm GL}}
\def\SL{{\rm SL}}

\def\Hom{{\rm Hom}}

\def\di!{\frac{\der^i}{i!}}
\def\dik!{\frac{\der^k_i}{k!}}

\def\N{\mathbb{N}}

\def\0{\overline{0}}
\def\1{\overline{1}}

\def\Ln1{\L_{n,\overline{1}}}

\def\a1{a_{\overline{1}}}

\def\S{\Sigma}

\def\vn1{\overrightarrow{n-1}}

\def\im{{\rm im}}

\def\mA{\mathbb{A}}

\def\Inn{{\rm Inn}}

\def\mS{\mathbb{S}}

\def\mT{\mathbb{T}}
\def\ind{{\rm ind}}

\begin{document}

\author{V. V. \  Bavula  
}

\title{${\rm K}_1(\mS_1)$ and the group of automorphisms of the  algebra $\mS_2$
 of one-sided inverses of a polynomial algebra in two variables }

\maketitle

\begin{abstract}

Explicit generators are found for the group $G_2$ of automorphisms
of the algebra $\mS_2$ of one-sided inverses of a polynomial
algebra in two variables over a field. 
Moreover, it is proved that
$$ G_2\simeq S_2\ltimes \mT^2\ltimes \Z\ltimes ((K^*\ltimes
E_\infty (\mS_1))\boxtimes_{\GL_\infty (K)}(K^*\ltimes E_\infty
(\mS_1)))$$ where $S_2$ is the symmetric group, $\mT^2$ is the
$2$-dimensional algebraic torus, $E_\infty (\mS_1)$ is the
subgroup of $\GL_\infty (\mS_1)$ generated by the elementary
matrices. In the proof, we use and prove several results on the
index of operators, and the final argument in the proof is  the
fact that ${\rm K}_1 (\mS_1) \simeq K^*$ proved in the paper. The
algebras $\mS_1$ and $\mS_2$ are noncommutative, non-Noetherian,
and  not domains. The group of units of the algebra $\mS_2$ is
found (it is huge).




\end{abstract}


\section{Introduction}
Throughout, ring means an associative ring with $1$; module means
a left module;
 $\N :=\{0, 1, \ldots \}$ is the set of natural numbers; $K$ is a
field 
and  $K^*$ is its group of units;
$P_n:= K[x_1, \ldots , x_n]$ is a polynomial algebra over $K$;
$\der_1:=\frac{\der}{\der x_1}, \ldots , \der_n:=\frac{\der}{\der
x_n}$ are the partial derivatives ($K$-linear derivations) of
$P_n$; $\End_K(P_n)$ is the algebra of all $K$-linear maps from
$P_n$ to $P_n$ and $\Aut_K(P_n)$ is its group of units (i.e. the
group of all the invertible linear maps from $P_n$ to $P_n$); the
subalgebra  $A_n:= K \langle x_1, \ldots , x_n , \der_1, \ldots ,
\der_n\rangle$ of $\End_K(P_n)$ is called the $n$'th {\em Weyl}
algebra.

$\noindent $

{\it Definition}, \cite{shrekalg}. The 
{\em algebra} $\mathbb{S}_n$ {\em of one-sided inverses} of $P_n$
is an algebra generated over a field $K$ 
 by $2n$ elements $x_1, \ldots , x_n, y_n, \ldots , y_n$ that
satisfy the defining relations:
$$ y_1x_1=\cdots = y_nx_n=1 , \;\; [x_i, y_j]=[x_i, x_j]= [y_i,y_j]=0
\;\; {\rm for\; all}\; i\neq j,$$ where $[a,b]:= ab-ba$, the
commutator of elements $a$ and $b$.

$\noindent $

By the very definition, the algebra $\mS_n$ is obtained from the
polynomial algebra $P_n$ by adding commuting, left (but not
two-sided) inverses of its canonical generators. The algebra
$\mS_1$ is a well-known primitive algebra \cite{Jacobson-StrRing},
p. 35, Example 2. Over the field
 $\mathbb{C}$ of complex numbers, the completion of the algebra
 $\mS_1$ is the {\em Toeplitz algebra} which is the
 $\mathbb{C}^*$-algebra generated by a unilateral shift on the
 Hilbert space $l^2(\N )$ (note that $y_1=x_1^*$). The Toeplitz
 algebra is the universal $\mathbb{C}^*$-algebra generated by a
 proper isometry.

$\noindent $

{\it Example}, \cite{shrekalg}. Consider a vector space $V=
\bigoplus_{i\in \N}Ke_i$ and two shift operators on $V$, $X:
e_i\mapsto e_{i+1}$ and $Y:e_i\mapsto e_{i-1}$ for all $i\geq 0$
where $e_{-1}:=0$. The subalgebra of $\End_K(V)$ generated by the
operators $X$ and $Y$ is isomorphic to the algebra $\mS_1$
$(X\mapsto x$, $Y\mapsto y)$. By taking the $n$'th tensor power
$V^{\t n }=\bigoplus_{\alpha \in \N^n}Ke_\alpha$ of $V$ we see
that the algebra $\mS_n\simeq \mS_1^{\t n}$ is isomorphic to the
subalgebra of $\End_K(V^{\t n })$ generated by the $2n$ shifts
$X_1, Y_1, \ldots , X_n, Y_n$ that act in different directions.

$\noindent $

 Let $G_n:=\Aut_{K-{\rm
alg}}(\mS_n)$ be the group of automorphisms of the algebra
$\mS_n$, and $\mS_n^*$ be the group of units of the algebra
$\mS_n$.

\begin{theorem}\label{Int24Apr9}
\begin{enumerate}
\item {\rm \cite{shrekaut}} $\; G_n=S_n\ltimes \mT^n\ltimes \Inn
(\mS_n)$.\item {\rm \cite{shrekaut}} $\; G_1=\mT^1\ltimes
\GL_\infty (K)$.
\end{enumerate}
\end{theorem}
where  $S_n=\{ s \in S_n\, | \, s (x_i) = x_{s(i)}, s(y_i)=
y_{s(i)}\}$ is the symmetric group, $\mT^n:=\{ t_\l \, | \, t_\l
(x_i) = \l _ix_i, t_\l (y_i) = \l_i^{-1}y_i, \l =(\l_i)\in
K^{*n}\}$ is the $n$-dimensional algebraic torus, $\Inn (\mS_n)$
is the group of inner automorphisms of the algebra $\mS_n$ (which
is huge), and $\GL_\infty (K)$ is the group of all the invertible
infinite dimensional matrices of the type $1+M_\infty (K)$ where
the algebra (without 1) of infinite dimensional matrices $M_\infty
(K) :=\varinjlim M_d(K)=\bigcup_{d\geq 1}M_d(K)$ is the injective
limit of matrix algebras. A semi-direct product $H_1\ltimes
H_2\ltimes \cdots \ltimes H_m$ of several groups means that
$H_1\ltimes (H_2\ltimes ( \cdots \ltimes (H_{m-1}\ltimes
H_m)\cdots )$.

The results of the papers \cite{Bav-Jacalg, shrekalg,  shrekaut,
jacaut} and the present paper show that (when ignoring
non-Noetherian property) the algebra $\mS_n$ belongs to the family
of algebras like the $n$'th Weyl algebra $A_n$,  the polynomial
algebra $P_{2n}$ and the Jacobian algebra $\mA_n$ (see
\cite{Bav-Jacalg, jacaut}). The structure of the group
$G_1=\mT^1\ltimes \GL_\infty (K)$ is  another confirmation of
`similarity' of the algebras $P_2$, $A_1$, and $\mS_1$. The groups
of automorphisms of the polynomial algebra $P_2$ and the Weyl
algebra $A_1$ (when ${\rm char}(K)=0$) were found by Jung
\cite{jung}, Van der Kulk \cite{kulk}, and Dixmier \cite{Dix}
respectively. These two groups have almost identical structure,
they are `infinite $GL$-groups' in the sense that they are
generated by the algebraic torus $\mT^1$ and by the obvious
automorphisms: $x\mapsto x+\l y^i$, $y\mapsto y$; $x\mapsto x$,
$y\mapsto y+ \l x^i$, where $i\in \N$ and $\l \in K$; which are
sort of `elementary infinite dimensional matrices' (i.e. `infinite
dimensional transvections`). The same picture as for the group
$G_1$. In prime characteristic, the group of automorphism of the
Weyl algebra $A_1$ was found by Makar-Limanov \cite{Mak-LimBSMF84}
(see also Bavula \cite{A1rescen}  for a different approach and for
further developments).

\begin{theorem}\label{aInt24Apr9}
\begin{enumerate}
\item $\mS_n^* = K^*\times (1+\ga_n)^*$ where the ideal $\ga_n$ of
the algebra $\mS_n$ is the sum of all the height 1 prime ideals of
the algebra $\mS_n$.  \item  The centre of the group $\mS_n^*$ is
$K^*$, and the centre of the group $(1+\ga_n)^*$ is $\{ 1 \}$.
\item The map $(1+\ga_n)^*\ra \Inn (\mS_n)$, $u\mapsto \o_u$, is a
group  isomorphism ($\o_u(a):=uau^{-1}$ for $a\in \mS_n$).
\end{enumerate}
\end{theorem}

The proof of this theorem is given at the end of Section
\ref{KTG1S} (another proof via the Jacobian algebras is given in
\cite{jacaut}).


To save on notation, we identify the groups $(1+\ga_n)^*$ and
$\Inn (\mS_n)$ via $u\mapsto \o_u$. Clearly, $\mS_2=\mS_1(1)\t
\mS_1(2)$ where $\mS_1(i):=K\langle x_i, y_i\rangle \simeq \mS_1$.
The algebra $\mS_2$ has only two height one prime ideals $\gp_1$
and $\gp_2$. Let $F_2:=\gp_1\cap \gp_2$.  The aim of the paper is
to find generators for the group $G_2$ (see the end of the paper)
and to prove the following theorem.
\begin{itemize}
\item (Theorem \ref{23Apr9})
\begin{enumerate}
 \item $ G_2= S_2\ltimes \mT^2\ltimes
\Theta \ltimes ((U_1(K)\ltimes E_\infty
(\mS_1(2)))\boxtimes_{(1+F_2)^*}(U_2(K)\ltimes E_\infty
(\mS_1(1)))$.
 \item $ G_2\simeq S_2\ltimes \mT^2\ltimes
\Z\ltimes ((K^*\ltimes E_\infty (\mS_1))\boxtimes_{\GL_\infty
(K)}(K^*\ltimes E_\infty (\mS_1)))$.
\end{enumerate}
\end{itemize}
For a ring $R$, let $E_\infty (R)$ be the subgroup of $\GL_\infty
(R)$ generated by  all the elementary matrices $1+rE_{ij}$ (where
$i,j\in \N$ with   $i\neq j$, and $r\in R$),  and let $U(R):= \{
uE_{00}+1-E_{00}\, | \, u\in R^*\}\simeq R^*$ where $R^*$ is  the
group of units of the ring $R$. The group $E_\infty (R)$ is a
normal subgroup of $\GL_\infty (R)$. The group $\Theta \simeq \Z$
is generated by a single element $\th$ (see Section \ref{KTG1S}).

 If a group $G$ is equal to the product $AB:=\{ ab\, | \, a\in A,
b\in B\}$ of its {\em normal} subgroups $A$ and $B$ then we write
$G=A \boxtimes B= A\boxtimes_{A\cap B}B$. So, each element $g\in
G$ is a product $ab$ for some elements $a\in A$ and $b\in B$; and
$ab=a'b'$ (where $a'\in A$ and $b'\in B$) iff $a'=ac$ and
$b'=c^{-1}b$ for some element $c\in A\cap B$. Clearly, $A\boxtimes
B = B\boxtimes A$.

At the final stage of the proof of  Theorem  \ref{23Apr9} we use
the fact that
\begin{itemize}
\item (Theorem \ref{d24Apr9}) ${\rm K}_1(\mS_1)\simeq K^*$.
\end{itemize}
The group of units $\mS_2^*$ of the algebra $\mS_2$ is found.

\begin{itemize}
\item (Corollary \ref{a25Apr9})  $ \mS_2^*= K^*\times \Theta
\ltimes ((U_1(K)\ltimes E_\infty
(\mS_1(2))\boxtimes_{(1+F_2)^*}(U_2(K)\ltimes E_\infty
(\mS_1(1)))$.
\end{itemize}

{\bf The structure of the proof of Theorem \ref{23Apr9}}. The
$\mS_2$-module $P_2$ is faithful, and so $\mS_2\subseteq
\End_K(P_2)$. By Theorem \ref{Int24Apr9} and Theorem
\ref{aInt24Apr9}, the question of finding the group $G_2=
S_2\ltimes \mT^2\ltimes  (1+\ga_2)^*$ is equivalent to the
question of finding the group $(1+\ga_2)^*$ or $\mS_2^* =
K^*\times (1+\ga_2)^*$. Difficulty in finding the group $\mS_2^*$
stems from two facts: (i) $\mS_2^*\subsetneqq \mS_2\cap
\Aut_K(P_2)$, i.e. there are non-units of the algebra $\mS_2$ that
are invertible linear maps in $P_2$; and (ii) some units of the
algebra $\mS_2$ are product of  {\em non-units} with {\em
non-zero} indices (each unit has zero index). To eliminate (ii)
the group $\Theta$ is introduced, and it is proved that
$(1+\ga_2)^* = \Theta \ltimes \CK$ and for the normal subgroup
$\CK $ of $(1+\ga_2)^*$ the situation (ii) does not occur. The
group $\CK$ is the common kernel of group epimorphisms $\ind_i :
(1+\ga_2)^*\ra \Z$ (see (\ref{indi})) where $i=1,2$. In order to
construct the maps $\ind_i$ and  to prove that they are well
defined group homomorphisms we need several results on the index
of operators which are collected at the beginning of Section
\ref{KTG1S}. Some of these are new (Theorem \ref{b23Apr9} and
Corollary \ref{b21Apr9}). Briefly, using indices of operators is
the main tool in finding the group $\mS_2^*$ and to prove that
${\rm K}_1(\mS_1)\simeq K^*$. Using indices  and the fact that
$(1+F_2)^* = (1+F_2)\cap \Aut_K(P_2)$ we show that $\CK =
(1+\gp_1)^*\boxtimes_{(1+F_2)^*}(1+\gp_2)^*$ (Proposition
\ref{b24Apr9}). Then using indices, we prove that $(1+\gp_i)^* =
U_i(K)\ltimes E_\infty (\mS_1(i+1))$ (Proposition \ref{c24Apr9}).
 This fact is equivalent to the fact  that ${\rm K}_1(\mS_1)
\simeq K^*$ (Theorem \ref{d24Apr9}).


\section{The groups $G_2$ and ${\rm K}_1(\mS_1)$}\label{KTG1S}

In this section, the groups $G_2$, $\mS_2^*$, and  and $K_1(\mS_1)
$ are found (Theorems \ref{23Apr9} and \ref{d24Apr9}, Corollary
\ref{a25Apr9}). The proofs are constructive.

We mentioned already in the Introduction that the key idea in
finding the group $G_2$ is to use indices of operators. That is
why we start this section with collecting known results on indices
and prove new ones. These results are used in all the proofs that
follow.

{\bf The index $\ind$ of linear maps and its properties}. Let $\CC
$ be the family of all $K$-linear maps with finite dimensional
kernel and cokernel, i.e. $\CC$ is the family of {\bf Fredholm}
linear maps/operators. For vector spaces $V$ and $U$, let $\CC
(V,U)$ be the set of all the linear maps from $V$ to $U$ with
finite dimensional kernel and cokernel. So, $\CC =\bigcup_{V,U}\CC
(V,U)$ is the disjoint union.

$\noindent $

{\it Definition}. For a linear map $\v \in \CC$, the integer $
\ind (\v ) := \dim \, \ker (\v ) - \dim \, \coker (\v )$ is called
the {\bf index} of the map $\v$.

$\noindent $

For vector spaces $V$ and $U$, let $\CC (V,U)_i:= \{ \v \in \CC
(V,U)\, | \, \ind (\v ) = i\}$. Then $\CC (V,U)=\bigcup_{i\in
\Z}\CC (V,U)_i$ is the disjoint union, and the family $\CC$ is the
disjoint union $\bigcup_{i\in \Z}\CC_i$ where $\CC_i :=\{ \v \in
\CC \, | \, \ind (\v ) = i\}$. When $V=U$, we write $\CC (V):= \CC
(V,V)$ and $\CC (V)_i:= \CC (V,V)_i$.

{\it Example}. Note that $\mS_1\subset \End_K(P_1)$
($x*x^i=x^{i+1}$, $y*x^{i+1}=x^i$, $i\in \N$, and $y*1=0$). The
map $x^i$ acting on the polynomial algebra $P_1$ is an injection
with $P_1=(\bigoplus_{j=0}^{i-1}Kx^j)\bigoplus \im (x^i)$; and the
map $y^i$ acting on $P_1$ is a surjection with $\ker (y^i)
=\bigoplus_{j=0}^{i-1}Kx^j$, and so  
\begin{equation}\label{indxy}
\ind (x^i)= -i\;\; {\rm and }\;\; \ind (y^i)= i, \;\; i\geq 1.
\end{equation}

Lemma \ref{b8Feb9} shows that $\CC$ is a multiplicative semigroup
with zero element (if the composition of two elements of $C$ is
undefined we set their product to be zero).

\begin{lemma}\label{b8Feb9}
Let $\psi : M\ra N $ and $\v : N\ra L$ be $K$-linear maps. If two
of the following three maps: $\psi$, $\v$,  and $\v \psi$, belong
to the set $\CC$ then so does the third; and in this case, $$ \ind
(\v\psi ) = \ind (\v ) + \ind (\psi ).$$
\end{lemma}
 By Lemma \ref{b8Feb9}, $\CC (N,L)_i\CC (M,N)_j\subseteq \CC
 (M,L)_{i+j}$ for all $i,j\in \Z$.
\begin{lemma}\label{a8Feb9}
Let
$$
\xymatrix{0\ar[r] & V_1\ar[r]\ar[d]^{\v_1}  & V_2 \ar[r]\ar[d]^{\v_2} & V_3 \ar[r]\ar[d]^{\v_3 } & 0 \\
0\ar[r] & U_1\ar[r]  & U_2\ar[r] & U_3 \ar[r] & 0 }
$$
be a commutative diagram of $K$-linear maps with exact rows.
Suppose that $\v_1, \v_2, \v_3\in \CC$. Then
$$ \ind (\v_2) = \ind (\v_1)+\ind (\v_3).$$
\end{lemma}

Let $V$ and $U$ be vector spaces. Define $\CI (V,U) :=\{ \v \in
\Hom_K(V,U)\, | \, \dim \, \im (\v ) <\infty \}$, and when $V=U$
we write $\CI (V):= \CI (V,V)$.

\begin{theorem}\label{b23Apr9}
Let $V$ and $U$ be vector spaces. Then $\CC (V,U)_i+\CI (V, U) =
\CC (V,U)_i$ for all $i\in \Z$.
\end{theorem}

{\it Proof}. The theorem is obvious if the vector spaces $V$ and
$U$ are finite dimensional since in this case the index of each
linear map from $V$ to $U$ is equal to $\dim (V) - \dim (U)$.  We
deduce the general case from this one. Let $u\in \CC (V,U)_i$ and
$f\in \CI (V,U)$. Using the fact that the kernel $\ker (f)$ has
finite codimension in $V$, i.e. $\dim (V/\ker (f))<\infty $ (since
 $V/\ker (f)\simeq \im (f)$), we can easily find subspaces $V_1,
 V_2\subseteq V$ and $W,U_1, U_2\subseteq U$ such that $\dim
 (V_1)<\infty $, $\dim (U_1)<\infty$, $\dim (W)<\infty$,
 $$V=\ker (u)\bigoplus
 V_1\bigoplus V_2, \; U= W \bigoplus U_1\bigoplus U_2, \;
 u|_{V_1}: V_1\simeq U_1, \; u|_{V_2}: V_2\simeq U_2,\;  f(V_2)
 =0$$
 and $f(\ker (u) \bigoplus V_1)\subseteq W \bigoplus U_1$. Note
 that $\im (u) = U_1\bigoplus U_2$ and $U/\im (u) \simeq W$.
 Consider the restrictions, say $u'$ and $f'$, of the maps $u$ and
 $f$ to the finite dimensional subspace $\ker (u) \bigoplus V_1$ of $V$,
 i.e.
\begin{equation}\label{u1f1}
u', f': \ker (u) \bigoplus V_1\ra W \bigoplus U_1.
\end{equation}
Then it is obvious that $\ind (u') = \ind (u)$ and $\ind (u'+f') =
\ind (u+f)$ (since $u+f|_{V_2}= u|_{V_2}:V_2\simeq U_2$). On the
other hand, $\ind (u'+f') = \ind (u')$ since the vector spaces in
(\ref{u1f1}) are finite dimensional. Therefore, $\ind (u+f) = \ind
(u)$. $\Box $

\begin{lemma}\label{a21Apr9}
Let $V$ and $V'$ be vector spaces, and  $\v :V\ra V'$ be a linear
map such that the vector spaces $\ker (\v )$ and $\coker (\v )$
 are isomorphic. Fix subspaces $U\subseteq V$ and $W\subseteq V'$
 such that $V=\ker (\v ) \bigoplus U$ and $V'= W\bigoplus \im (\v
 )$  and fix an isomorphism $f: \ker (\v ) \ra W$ (this is
 possible since $\ker (\v ) \simeq \coker (\v ) \simeq W$) and
 extend it to a linear map $f: V\ra V'$ by setting $f(U)=0$. Then
 the map $\v +f:V\ra V'$ is an isomorphism.
\end{lemma}

{\it Proof}. The map $\v +f$ is a surjection since $(\v +f)(V) =
(\v +f) (\ker (\v ) +U) = W+\im (\v ) = V'$. The map $\v +f$ is an
injection: if $v\in \ker (\v +f)$ then $\v (v) = f(-v) \in W\cap
\im (\v ) =0$, and so $v\in \ker (\v ) \cap \ker (f) = \ker (\v )
\cap U=0$. Therefore, the map $\v +f$ is an isomorphism. $\Box $

\begin{lemma}\label{b21Apr9}
Let $V$ and $V'$ be vector spaces,  $\v \in \CC ( V, V')_i$ for
some $i\in \Z$, $V=\ker (\v ) \bigoplus U$ and $V'=W\bigoplus \im
(\v )$
 for some subspaces $U\subseteq V$ and $W\subseteq V'$.
\begin{enumerate}
\item If $\dim \, \ker (\v ) \leq \dim \, \coker (\v )$    then
fix an injective linear map $f: \ker (\v ) \ra W$ and extend it to
a linear map $f: V\ra V'$ by setting $f(U)=0$. Then the map $\v
+f$ is an injection that belongs to  $\CC ( V, V')_i$.  \item If
$\dim \, \ker (\v ) \geq \dim \, \coker (\v )$ then fix a
surjective linear map $f: \ker (\v ) \ra W$ and extend it to a
linear map $f: V\ra V'$ by setting $f(U)=0$. Then the map $\v +f$
is a surjection that belongs to  $\CC ( V, V')_i$.
\end{enumerate}
\end{lemma}

{\it Proof}. 1. An arbitrary element $v\in V= \ker (\v ) \bigoplus
U$ is a unique sum $k+u$ where $k\in \ker (\v )$ and $u\in U$. If
$v\in \ker (\v+f)$ then $0=(\v +f) (k+u) = f(k) +\v (u)$ and so
$f(k)=0$ and $\v (u)=0$ (since $V'=W\bigoplus \im (\v )$, $f(k)\in
W$ and $\v (u) \in \im (\v )$) hence $k=0$ ($f$ is an injection)
and $u=0$ ($\v|_{U}:U\simeq \im (\v )$). Therefore, $v=0$. This
means that the map $\v +f$ is an injection. Now, $\im (\v +f) =
f(\ker (\v ))\bigoplus \im (\v )$, and so
$$ \ind (\v +f) = -\dim \, \coker (\v +f) = -\dim (W/\ker (\v )) =
\dim \, \ker (\v ) -\dim (W) = \ind (\v ).$$

2. The map $\v +f$ is a surjection since $(\v +f) (V) = (\v
+f)(\ker (\v ) +U) = f(\ker (\v ))+\im (\v ) = W+\im (\v ) = V'$.
Suppose that an element $v= k+u\in V$ (where $k\in \ker (\v )$ and
$u\in U$) belongs to the kernel of the map $\v +f$. Then $0=(\v
+f) (v) = f(k) +\v (u)$ and so $f(k) =0$ and $\v (u) =0$ (since
$V'=W\bigoplus \im (\v )$, $f(k) \in W$ and $\v (u) \in \im (U)$),
hence $k\in \ker (f)$ and $u=0$ since ($\v |_U: U\simeq \im (\v
)$). Now,
$$ \ind (\v +f) = \dim \, \ker (f) =   \dim \, \ker (\v)  -\dim (W)= \ind (\v ).\;\;\; \Box$$

$\noindent $

{\bf The algebras $\mS_1$ and $\mS_2$}. We collect some results
without proofs on the algebras $\mS_1$ and  $\mS_2$ from
\cite{shrekalg} that will be used in this paper, their proofs can
be found in \cite{shrekalg}. Clearly, $\mathbb{S}_2=\mS_1(1)\t
\mS_1(2)\simeq \mathbb{S}_1^{\t 2}$ where $\mS_1(i):=K\langle
x_i,y_i \, | \, y_ix_i=1\rangle \simeq \mS_1$ and
$\mS_2=\bigoplus_{\alpha , \beta \in \N^2} Kx^\alpha y^\beta$
where $x^\alpha := x_1^{\alpha_1} x_2^{\alpha_2}$, $\alpha =
(\alpha_1, \alpha_2)$, $y^\beta := y_1^{\beta_1} y_2^{\beta_2}$,
$\beta = (\beta_1, \beta_2)$. In particular, the algebra $\mS_2$
contains two polynomial subalgebras $P_2$ and $Q_2:=K[y_1, y_2]$
and is equal, as a vector space,  to their tensor product $P_2\t
Q_2$.

When $n=1$, we usually drop the subscript `1' if this does not
lead to confusion.  So, $\mS_1= K\langle x,y\, | \,
yx=1\rangle=\bigoplus_{i,j\geq 0}Kx^iy^j$. For each natural number
$d\geq 1$, let $M_d(K):=\bigoplus_{i,j=0}^{d-1}KE_{ij}$ be the
algebra of $d$-dimensional matrices where $\{ E_{ij}\}$ are the
matrix units,  $M_\infty (K) := \varinjlim
M_d(K)=\bigoplus_{i,j\in \N}KE_{ij}$ be the algebra (without 1) of
infinite dimensional matrices, and $\GL_\infty (K)$ be the  group
of units of the monoid $1+M_\infty (K)$. The algebra $\mS_1$
contains the ideal $F:=\bigoplus_{i,j\in \N}KE_{ij}$, where
\begin{equation}\label{Eijc}
E_{ij}:= x^iy^j-x^{i+1}y^{j+1}, \;\; i,j\geq 0.
\end{equation}
For all natural numbers $i$, $j$, $k$, and $l$,
$E_{ij}E_{kl}=\d_{jk}E_{il}$ where $\d_{jk}$ is the Kronecker
delta function.  The ideal $F$ is an algebra (without 1)
isomorphic to the algebra $M_\infty (K)$ via $E_{ij}\mapsto
E_{ij}$.
\begin{equation}\label{mS1d}
\mS_1= K\oplus xK[x]\oplus yK[y]\oplus F,
\end{equation}
the direct sum of vector spaces. Then 
\begin{equation}\label{mS1d1}
\mS_1/F\simeq K[x,x^{-1}]=:L_1, \;\; x\mapsto x, \;\; y \mapsto
x^{-1},
\end{equation}
since $yx=1$, $xy=1-E_{00}$ and $E_{00}\in F$. The algebra $\mS_2
= \mS_1(1)\t \mS_1(2)$ contains the ideal
$$F_2:= F(1) \t F(2)=\bigoplus_{\alpha , \beta \in
\N^2}KE_{\alpha \beta}, \;\; {\rm where}\;\; E_{\alpha
\beta}:=\prod_{i=1}^2 E_{\alpha_i \beta_i}(i),$$ where $F(i)$ is
the ideal $F$ of the algebra $\mS_1(i)$ and
$E_{\alpha_i\beta_i}(i)$ are its matrix units as defined in
(\ref{Eijc}). Note that $E_{\alpha \beta}E_{\g \rho}=\d_{\beta \g
}E_{\alpha \rho}$ for all elements $\alpha, \beta , \g , \rho \in
\N^2$ where $\d_{\beta
 \g }$ is the Kronecker delta function, and so $(1+F_2)^*\simeq \GL_\infty (K)$.

 The algebra $\mS_2$ contains only two height one prime ideals
 $\gp_1:= F(1) \t \mS_1(2)$ and $\gp_2:= \mS_1(1)\t F(2)$.
 Clearly, $F_2=\gp_1\cap \gp_2= \gp_1\gp_2$, $\mS_2/ \gp_1\simeq
 K[x_1, x_1^{-1}]\t \mS_1(2)$ and $\mS_2/ \gp_2\simeq \mS_1(1) \t
 K[x_2,x_2^{-1}]$. The ideal $\ga_2:= \gp_1+\gp_2$ plays an
 important role in this paper, $\mS_2/\ga_2\simeq K[x_1, x_1^{-1}]\t
 K[x_2,x_2^{-1}]$.
\begin{proposition}\label{a19Dec8}
{\rm \cite{shrekalg}} The polynomial algebra $P_n$
 is the only (up to isomorphism)  faithful, simple $\mS_n$-module.
\end{proposition}

In more detail, ${}_{\mS_n}P_n\simeq \mS_n / (\sum_{i=0}^n \mS_n
y_i) =\bigoplus_{\alpha \in \N^n} Kx^\alpha \overline{1}$,
$\overline{1}:= 1+\sum_{i=1}^n \mS_ny_i$; and the action of the
canonical generators of the algebra $\mS_n$ on the polynomial
algebra $P_n$ is given by the rule:
$$ x_i*x^\alpha = x^{\alpha + e_i}, \;\; y_i*x^\alpha = \begin{cases}
x^{\alpha - e_i}& \text{if } \; \alpha_i>0,\\
0& \text{if }\; \alpha_i=0,\\
\end{cases}  \;\; {\rm and }\;\; E_{\beta \g}*x^\alpha = \d_{\g
\alpha} x^\beta,
$$
where the set of elements $e_1:= (1,0,\ldots , 0),  \ldots ,
e_n:=(0, \ldots , 0,1)$ is the canonical basis for the free
$\Z$-module $\Z^n=\bigoplus_{i=1}^n \Z e_i$.  We identify the
algebra $\mS_n$ with its image in the algebra $\End_K(P_n)$ of all
the $K$-linear maps from the vector space $P_n$ to itself, i.e.
$\mS_n \subset \End_K(P_n)$.

\begin{corollary}\label{a23Apr9}
\begin{enumerate}
\item $1+F_2\subseteq \CC (P_2)_0$. \item $\mS_2^* +F_2\subseteq
\CC (P_2)_0$.
\end{enumerate}
\end{corollary}

{\it Proof}. Both statements follows from Theorem \ref{b23Apr9}:
$\mS_2^*\in \CC (P_2)_0$ and $F_2\in \CI (P_2)$, but we give short
independent proofs (that do not use Theorem \ref{b23Apr9}).

 1. Since $1+F_2\simeq 1+M_\infty (K)$, statement 1 is obvious.

2. Let $u\in \mS_2^*$ and $f\in F_2$. Then $u^{-1}f\in F_2$. By
statement 1, the element $1+u^{-1}f\in \CC (P_2)_0$. Since $u\in
\CC (P_2)_0$, we have $u+f= u(1+u^{-1}f)\in \CC (P_2)_0$, by Lemma
\ref{b8Feb9}.  $\Box $

$\noindent $

{\bf The subgroup $\Theta$ of $(1+\ga_2)^*$}. The element $\th :=
(1+(y_1-1)E_{00}(2))(1+E_{00}(1)(x_2-1))\in (1+\ga_2)^*$ is a unit
and 
\begin{equation}\label{th12el}
\th^{-1}= (1+E_{00}(1)(y_2-1))(1+(x_1-1)E_{00}(2))\in (1+\ga_2)^*.
\end{equation}
This is obvious since
$$\th *x^\alpha = \begin{cases}
x^\alpha& \text{if } \alpha_1>0, \alpha_2>0,\\
x_1^{\alpha_1-1}& \text{if } \alpha_1>0, \alpha_2=0,\\
x_2^{\alpha_2+1}& \text{if } \alpha_1=0, \alpha_2\geq 0, \\
\end{cases}
\;\;\; {\rm and}\;\;\; \th^{-1} *x^\alpha = \begin{cases}
x^\alpha& \text{if } \alpha_1>0, \alpha_2>0,\\
x_1^{\alpha_1+1}& \text{if } \alpha_1\geq 0, \alpha_2=0,\\
x_2^{\alpha_2-1}& \text{if } \alpha_1=0, \alpha_2> 0. \\
\end{cases}
$$
Let $\Theta$ be the subgroup of $(1+\ga_2)^*$ generated by the
element $\th$. Then $\Theta \simeq \Z$ since $\th^i*1=x_2^i$ for
all $i\geq 1$. It follows from
$$ (1+(y_1-1)E_{00}(2))*x^\alpha = \begin{cases}
x^\alpha& \text{if } \alpha_2>0,\\
x_1^{\alpha_1-1}& \text{if } \alpha_1>0, \alpha_2=0,\\
0& \text{if } \alpha_1=0, \alpha_2= 0, \\
\end{cases}$$
that the map $1+(y_1-1)E_{00}(2)\in \End_K(P_2)$ is a surjection
with kernel equal to $K$, and so 
\begin{equation}\label{indEy}
\ind (1+(y_1-1)E_{00}(2)) = 1.
\end{equation}
Similarly, it follows from
$$ (1+E_{00}(1)(x_2-1))*x^\alpha = \begin{cases}
x^\alpha& \text{if } \alpha_1>0,\\
x_2^{\alpha_2+1}& \text{if } \alpha_1=0,\\
\end{cases}$$
that the map $1+E_{00}(1)(x_2-1)\in \End_K(P_2)$ is an injection
with $P_2= K\bigoplus \im (1+E_{00}(1)(x_2-1))$, and so
\begin{equation}\label{indEy1}
\ind (1+E_{00}(1)(x_2-1)) = -1.
\end{equation}
We see that the unit $\th$ of the algebra $\mS_2$ is the product
of two {\em non-units} having nonzero indices  the sum of which is
equal to zero since  $\ind (\th )=0$. Lemma \ref{a24Apr9} shows
that this is a general phenomenon, and so the group $(1+\ga_2)^*$
is  a sophisticated group in the sense that in construction of
units non-units are involved.
\begin{lemma}\label{a24Apr9}
Let $u=1+a_1+a_2\in (1+\ga_2)^*$ where $a_i\in \gp_i$. Then
\begin{enumerate}
\item $1+a_1, 1+a_2\in \CC (P_2)$ and $\ind (1+a_1) +\ind
(1+a_2)=0$. \item If $u=1+a_1'+a_2'$ where $a_i'\in \gp_i$ then
$\ind (1+a_1) = \ind (1+a_1')$ and  $\ind (1+a_2) = \ind
(1+a_2')$.
\end{enumerate}
\end{lemma}

{\it Proof}. 1. Since $a_i\in \gp_i$, we have $a_1a_2, a_2a_1\in
F_2$. By Corollary \ref{a23Apr9}.(2), $u+a_1a_2, u+a_2a_1\in \CC
(P_2)_0$. Then, it follows from the equalities $u+a_1a_2=
(1+a_1)(1+a_2)$ and $u+a_2a_1= (1+a_2)(1+a_1)$, that
\begin{eqnarray*}
 \im (1+a_1) &\supseteq & \im ( u+a_1a_2), \; \;\; \ker (1+a_1) \subseteq  \ker ( u+a_2a_1),\\
 \im (1+a_2) &\supseteq & \im ( u+a_2a_1), \; \;\; \ker (1+a_2) \subseteq  \ker ( u+a_1a_2).\\
\end{eqnarray*}
This means that $1+a_1, 1+a_2\in \CC (P_2)$. By Corollary
\ref{a23Apr9}.(2) and Lemma \ref{b8Feb9}, $$0=\ind (u+a_1a_2) =
\ind (1+a_1)(1+a_2) = \ind (1+a_1) +\ind (1+a_2).$$

2. It is obvious that $a'= a_1+f$ and $a_2'= a_2-f$ for an element
$f\in \gp_1\cap \gp_2= F_2$. Since $F_2\subseteq \CI (P_2)$, we
see that $\ind (1+a_1') = \ind (1+a_1+f) = \ind (1+a_1)$ and $\ind
(1+a_2') = \ind (1+a_2-f) = \ind (1+a_2)$, by Theorem
\ref{b23Apr9}. $\Box $

$\noindent $

By Lemma \ref{a24Apr9}, for each number $i=1,2$, there is a
well-defined map,
\begin{equation}\label{indi}
\ind_i: (1+\ga_2)^*\ra \Z , \;\; u=1+a_1+a_2\mapsto \ind (1+a_i),
\end{equation}
which is a group homomorphism:
\begin{eqnarray*}
 \ind_i(uu') &= & \ind_i((1+a_1+a_2)(1+a_1'+a_2') )= \ind(1+a_i+a_i'+a_ia_i')     \\
 &=& \ind ((1+a_i)(1+a_i') )= \ind (1+a_i)+\ind(1+a_i')\\
 & =& \ind_i(u) +\ind_i(u').
\end{eqnarray*}
The restriction of the homomorphism $\ind_i$ to the subgroup
$\Theta$ is an isomorphism: $\ind_1: \Theta \ra \Z$, $\th \mapsto
-1$; $\ind_2: \Theta \ra \Z$, $\th \mapsto 1$. Therefore, the
homomorphisms $\ind_i$ are epimorphisms which have the {\em same}
kernel (Lemma \ref{a24Apr9}.(1)) which we denote by $\CK$. Then,
\begin{equation}\label{1aatCK}
(1+\ga_2)^* = \Theta \ltimes \CK .
\end{equation}

It is obvious that
$(1+\gp_1)^*\boxtimes_{(1+F_2)^*}(1+\gp_2)^*\subseteq \CK$ since
 $(1+\gp_1)^* \cap (1+\gp_2)^*=(1+F_2)^*$.

\begin{proposition}\label{b24Apr9}
$\CK =(1+\gp_1)^*\boxtimes_{(1+F_2)^*}(1+\gp_2)^*$.
\end{proposition}

{\it Proof}. It suffices to show that each element $u=1+a_1+a_2$
of the group $\CK$ is a product $u_1u_2$ for some elements $u_i\in
(1+\gp_i)^*$. Note that  $1+a_1\in \CC (P_2)_0$. Fix a subspace,
say $W$, of $P_2$ such that $P_2= \ker (1+a_1) \bigoplus W$ and
$W=\bigoplus_{\alpha \in I} Kx^\alpha$ where $I$ is a subset of
$\N^2$. By Lemma \ref{a21Apr9}, we can find an element $f_1\in
F_2$ (since $\dim \, \ker (1+a_1) <\infty $, $W$ has a monomial
basis, and $f_1(W)=0$) such that $u_1:= 1+a_1+f_1\in \Aut_K(P_2)$.
We claim that $u_1\in (1+\gp_1)^*$. It is a subtle point since not
all elements of the algebra $\mS_2$ that are invertible linear
maps in $P_2$ are invertible in $\mS_2$, i.e. $\mS_2^* \subsetneqq
\mS_2\cap \Aut_K(P_2)$ but $(1+F_2)^* = (1+F_2)\cap \Aut_K(P_2)$,
\cite{shrekaut}. The main idea in the proof of the claim is to use
this equality. Similarly, we can find an element $f_2\in F_2$ such
that $v:= 1+a_2+f_2\in \Aut_K(P_2)$. Then $u= u_1v+g_1$ and $u=
vu_1+g_2$ for some elements $g_i\in F_2$. Hence,
$$ u_1vu^{-1}=1-g_1u^{-1}\;\; {\rm and }\;\;\ u^{-1} vu_1= 1-u^{-1}
g_2, $$ and so $1-g_1u^{-1}, 1-u^{-1} g_2\in (1+F_2) \cap
\Aut_K(P_2) = (1+F_2)^*$. It follows that $u_1^{-1} = vu^{-1}
(1-g_1u^{-1})^{-1}\in (1+\gp_1)^*$ since
$$1\equiv  1-g_1u^{-1}
\equiv u_1vu^{-1} \equiv vu^{-1}\mod \gp_1.$$ This proves the
claim. Clearly, $u_2:= v+u_1^{-1}g_1\in 1+\gp_2$. Then, it follows
from the equality $ u = u_1v+g_1 = u_1(v+u_1^{-1}g_1)= u_1u_2$
that $u_2= u_1^{-1} u\in (1+\gp_2)^*$. This finishes the proof of
the proposition. $\Box $

$\noindent $

In order to save on notation, it is convenience to treat the set
of indices $\{ 1, 2\}$ as the group $\Z / 2\Z= \{ 1,2\}$ where
$1+1=2$ and $1+2=1$. For each number $i=1,2$, the group of units
of the monoid $1+\gp_i = 1+F(i) \t \mS_1(i+1)= 1+M_\infty
(\mS_1(i+1))$ is equal to $(1+\gp_i)^* = \GL_\infty (\mS_1(i+1))$.
It contains the semi-direct product $U_i(K)\ltimes E_\infty
(\mS_1(i+1))$ of its two subgroups, where
$$U_i(K):=\{ \l E_{00}(i) +1-E_{00}(i)\, | \, \l \in K^*\} \simeq
K^*$$
 and the group $E_\infty (\mS_1(i+1))$ is generated by all the
 elementary matrices $1+aE_{kl}(i)$ where $k\neq l$ and $a\in
 \mS_1(i+1)$. Note that the group $E_\infty (\mS_1(i+1))$ is a
 normal subgroup  of $\GL_\infty (\mS_1(i+1))$.

The set $F_2$ is an ideal of the algebra $K+\gp_i= K(1+\gp_i)$
which is a subalgebra  of the algebra $\mS_2$, and $(K+\gp_i) /
F_2= K(1+\gp_i/ F_2) \simeq K(1+M_\infty (L_{i+1}))$ where
$L_{i+1}:= K[x_{i+1}, x_{i+1}^{-1}]\simeq \mS_1(i+1)/F(i+1)$ is
the Laurent polynomial algebra. The algebra $L_{i+1}$ is a
Euclidian  domain, hence $\GL_\infty (L_{i+1}) = U(L_{i+1})\ltimes
E_\infty (L_{i+1})$ where
$$U(L_{i+1}):=\{ a
E_{00}(i) +1-E_{00}(i)\, | \, a \in L_{i+1}^*\} \simeq
L_{i+1}^*=K^*\times \{ x_{i+1}^m \, | \, m\in \Z\}$$ and $E_\infty
(L_{i+1})$ is the subgroup of $\GL_\infty (L_{i+1})$ generated by
all the elementary matrices.

The group of units of the algebra $(K+\gp_i) / F_2$ is equal to
$K^*\times \GL_\infty (L_{i+1})= K^* \times ( U(L_{i+1})\ltimes
E_\infty (L_{i+1}))$. The algebra epimorphism $\psi_i : K+\gp_i\ra
(K+\gp_i)/F_2$, $a\mapsto a+F_2$, induces the exact sequence of
groups, 
\begin{equation}\label{1piex}
1\ra (1+F_2)^*\ra (1+\gp_i)^*\stackrel{\psi_i }{\ra} \GL_\infty
(L_{i+1})=U(L_{i+1})\ltimes E_\infty (L_{i+1}),
\end{equation}
which yields the short exact sequence of groups, 
\begin{equation}\label{2piex}
1\ra (1+F_2)^*\ra U_i(K)\ltimes E_\infty (\mS_1(i+1))\ra
U(K)\ltimes E_\infty (L_{i+1})\ra 1
\end{equation}
since $(1+F_2)^* \subseteq E_\infty (\mS_1(i+1))$, by Proposition
\ref{10Apr10}.  Note that $U_i(K)\ltimes E_\infty
(\mS_1(i+1))\subseteq (1+\gp_i)^*$.  In fact, the equality holds.
\begin{proposition}\label{c24Apr9}
$(1+\gp_i)^* =U_i(K)\ltimes E_\infty (\mS_1(i+1))$.
\end{proposition}

{\it Proof}. In view of the exact sequences (\ref{1piex}) and
(\ref{2piex}), it suffices to show that the image of the map
$\psi_i$ in (\ref{1piex}) is equal to $U(K)\ltimes E_\infty
(L_{i+1}))$. Since
$$U(L_{i+1}) = U(K)\times \{E_{00}(i) x_{i+1}^m +1-E_{00}(i)  \, | \,
m\in \Z\},$$
 this is equivalent  to show that that if $\psi_i (u) =
E_{00}(i) x_{i+1}^m +1-E_{00}(i)$ for some element $u\in
(1+\gp_i)^*$ and an integer $m\in \Z$ then $m=0$. Let $u(m) :=
E_{00}(i) v_{i+1}(m)+1-E_{00}(i)$ where $v_{i+1}(m):=
\begin{cases}
x_{i+1}^m& \text{if } m\geq 0,\\
y_{i+1}^{|m|} & \text{if } m<0.\\
\end{cases} $ Then $u(m) \in 1+\gp_i$ and $\psi_i (u(m)) = \psi_i
(u)$. Hence, $u(m) = u+f_m$ for some element $f_m\in F_2$. Note
that $$u(m) =\begin{cases}
u(1)^m& \text{if } m\geq 0,\\
u(-1)^{|m|} & \text{if } m<0, \\
\end{cases} $$ and, by (\ref{indEy}) and (\ref{indEy1}), $\ind
(u(m))=-m$. By Corollary \ref{a23Apr9}.(2).
$$ 0=\ind (u) = \ind (u+f_m) = \ind (u(m)) = -m,$$
and so $m=0$, as required. $\Box $

$\noindent $

Proposition  \ref{c24Apr9} is equivalent to the next theorem.

\begin{theorem}\label{d24Apr9}
${\rm K}_1(\mS_1)\simeq U(K) \simeq K^*$ and $\GL_\infty
(\mS_1)=U(K)\ltimes E_\infty (\mS_1)$.
\end{theorem}

{\it Proof}. Recall that  ${\rm K}_1(\mS_1) := \GL_\infty (\mS_1)/
E_\infty (\mS_1)$ where $E_\infty (\mS_1)$ is the  subgroup of
$\GL_\infty (\mS_1)$ generated by the elementary matrices. The
group  $E_\infty (\mS_1)$ is a normal subgroup of $\GL_\infty
(\mS_1)$. It follows from $(1+\gp_i)^*\simeq \GL_\infty (\mS_1)$
and Proposition \ref{c24Apr9} that
$$ {\rm K}_1(\mS_1) \simeq U(K)\ltimes E_\infty (\mS_1)/ E_\infty
(\mS_1)\simeq U(K)\simeq K^*. \;\;\;  \Box $$

{\bf The determinant $\overline{\det}$}.  The algebra epimorphism
$\mS_1\ra \mS_1/F\simeq L_1$, $a\mapsto a+F$, yields the group
homomorphism $ \psi : \GL_\infty (\mS_1) \ra \GL_\infty (L_1)$. By
Proposition \ref{c24Apr9}, the image of the group homomorphism
$\det \circ \psi : \GL_\infty (\mS_1)
\stackrel{\psi}{\ra}\GL_\infty (L_1) \stackrel{\det}{\ra}L_1^*$ is
$K^*$. Therefore, there is a well determined group epimorphism:
\begin{equation}\label{bdet}
\overline{\det}:= \det \circ \psi : \GL_\infty (\mS_1)\ra K^*.
\end{equation}
By the very definition, $\overline{\det}(E_\infty (\mS_1))=1$ and
$\overline{\det}(\mu (\l ))=\l$ for all elements $\mu (\l )=\l
E_{00}+1-E_{00}\in U(K)$ where $\l \in K^* $.  Therefore, there is
the exact sequence of groups: 
\begin{equation}\label{1bdet}
1\ra E_\infty (\mS_1)\ra \GL_\infty
(\mS_1)\stackrel{\overline{\det}}{\ra}K^* \simeq {\rm
K}_1(\mS_1)\ra 1.
\end{equation}
\begin{corollary}\label{a11Apr10}
Each element $a$ of the group $\GL_\infty (\mS_1) = U(K)\ltimes
E_\infty (\mS_1)$ is a unique product $a= \mu (\l ) e$ where $\mu
(\l ) \in U(K)$ and $e\in E_\infty (\mS_1)$. Moreover, $\l =
\overline{\det}(a)$ and $e= \mu (-\overline{\det}(a)) a$.
\end{corollary}

 Recall that that $G_2=S_2\ltimes \mT^2\ltimes \Inn (\mS_2)$
 (Theorem \ref{Int24Apr9}.(1)) and  $(1+\ga_2)^*\simeq \Inn
(\mS_2)$, $u\lra \o_u$ (Theorem \ref{aInt24Apr9}.(3)). We {\em
identify} the groups $(1+\ga_2)^*$ and $\Inn (\mS_2)$ via $u\lra
\o_u$.

\begin{theorem}\label{23Apr9}
\begin{enumerate}
 \item $ G_2= S_2\ltimes \mT^2\ltimes
\Theta \ltimes ((U_1(K)\ltimes E_\infty
(\mS_1(2))\boxtimes_{(1+F_2)^*}(U_2(K)\ltimes E_\infty
(\mS_1(1)))$.
 \item $ G_2\simeq S_2\ltimes \mT^2\ltimes
\Z\ltimes ((K^*\ltimes E_\infty (\mS_1))\boxtimes_{\GL_\infty
(K)}(K^*\ltimes E_\infty (\mS_1)))$.
\end{enumerate}
\end{theorem}

{\it Proof}. By (\ref{1aatCK}), Proposition \ref{b24Apr9}, and
Proposition \ref{c24Apr9}, 
\begin{equation}\label{1a2}
(1+\ga_2)^* = \Theta \ltimes ((U_1(K)\ltimes E_\infty
(\mS_1(2))\boxtimes_{ (1+F_2)^*}(U_2(K)\ltimes E_\infty
(\mS_1(1))),
\end{equation}
and the statements follow since $(1+F_2)^* \simeq \GL_\infty (K)$.
$\Box $

\begin{corollary}\label{a25Apr9}
$\mS_2^*= K^*\times\Theta \ltimes ((U_1(K)\ltimes E_\infty
(\mS_1(2))\boxtimes_{(1+F_2)^*}(U_2(K)\ltimes E_\infty
(\mS_1(1)))$.
\end{corollary}

{\bf Generators for the group $G_2$}. Using Theorem \ref{23Apr9}
and (\ref{mS1d}), we can easily write down a set of generators for
the group $G_2$:
\begin{eqnarray*}
 s:& x_i\mapsto x_{i+1}, y_i\mapsto y_{i+1}, \;\; i=1,2;  \\
 t_{(\l , 1)}: & x_1\mapsto \l x_1, \;\; y_1\mapsto \l^{-1}y_1,
 \;\;x_2\mapsto x_2, \;\; y_2\mapsto y_2, \;\;
 \l \in K^*;
\end{eqnarray*}
$\o_\th$, $\o_{\l E_{00}(1)+1-E_{00}(1)}$, $\o_{\l
E_{ij}(1)E_{kl}(2) +1- E_{ij}(1)E_{kl}(2)}$,  $\o_{1+\mu
x_2^mE_{ij}(1)}$, $\o_{1+\mu y_2^mE_{ij}(1)}$, $\o_{1+\mu
E_{ij}(1)}$ and   where $\l \in K^*$, $\mu \in K$,  $ m,i,j,k,l\in
\N$, $m\geq 1$, $i\neq j$
 (note that $s\o_{1+\l x_1^mE_{ij}(2)}s^{-1}=
\o_{1+\l x_2^mE_{ij}(1)}$, etc).

$\noindent $

Each element $\s \in G_n = S_n \ltimes \mT^n \ltimes \Inn (\mS_n)$
is a unique product $st_\l \o_u$ where $s\in S_n$, $t_\l \in
\mT^n$ and $\o_u\in \Inn (\mS_n)$. In \cite{jacaut}, for each
element $\s \in G_n$, using the elements $\s (x_i), \s(y_i)\in
\mS_n$, $i=1, \ldots , n$ explicit algebraic formulae are found
for the components $s$, $t_\l$, and $\o_u$ of $\s$. So, the
automorphism $\s $ can be effectively (in finitely many steps)
decomposed into the product $st_\l \o_u$.

\begin{proposition}\label{10Apr10}
$(1+F_2)^*\subseteq E_\infty (\mS_1(i))$ for $i=1,2$.
\end{proposition}

{\it Proof}. Due to symmetry it suffices to show that
$(1+F_2)^*\subseteq E_\infty (\mS_1(2))$. Recall that the group
$(1+F_2)^* = (1+\sum_{\alpha, \beta \in \N^2}KE_{\alpha\beta})$ is
non-canonically isomorphic to the group $\GL_\infty (K)$. To see
this we have to choose a bijection $b:\N^2\ra \N$. Then the matrix
units $E_{\alpha \beta}$ can be seen as the usual matrix units
$E_{b(\alpha )b(\beta )}$ and so $(1+F_2)^* \simeq \GL_\infty
(K)$. This isomorphism depends on the choice of the bijection.
Since $(1+F_2)^* \simeq \GL_\infty (K)$, the group $(1+F_2)^*$ is
generated by the elements $a=1+\l E_{ij}(1)E_{kl}(2)$ where $\l
\in K$ and $(i,k) \neq (j,l)$, and $b=1+\l E_{00}(1)E_{00}(2)$
where $\l \in K\backslash \{ -1\}$. It suffices to show that these
generators belong to the group $E_\infty (\mS_1(2))$.

 First, let us show that $a\in E_\infty (\mS_1(2))$. If $i\neq j$
 then obviously the inclusion holds. If $i=j$, i.e. $a=1+\l
 E_{ii}(1)E_{kl}(2)$, then necessarily $k\neq l$ since $(i,k)\neq
 (i,l)$. For an element $g$ and $h$ of a group, $[g,h]:=
 ghg^{-1}h^{-1}$ is their {\em group commutator}. For any natural
 number $l$ such that $l\neq i$, the elements
 $1+E_{il}(1)E_{kk}(2)$ and $1+\l E_{li}(1)E_{kl}(2)$ belong to
 the group $E_\infty (\mS_1(2))$. Then so does their commutator
\begin{equation}\label{comEi}
[1+E_{il}(1)E_{kk}(2), 1+\l E_{li}(1)E_{kl}(2)]=1+\l
E_{ii}(1)E_{kl}(2).
\end{equation}
Therefore, all the generators $a$ belong to the group $E_\infty
(\mS_1(2))$.

It remains to prove that $b\in E_\infty (\mS_1(2))$. In the
$2\times 2$ matrix ring $M_2(\mS_1(2))$ with entries in the
algebra $\mS_1(2)$ we have the equality, for all scalars $\l \in
K\backslash \{ -1\}$: 
\begin{equation}\label{comEi1}
\Bigl(
\begin{matrix} 1  & 0\\ -\frac{y_2}{1+\l } & 0
\end{matrix}\Bigr)\,
\Bigl(
\begin{matrix} 1 & \l x_2 \\  0& 1
\end{matrix}\Bigr)\,
\Bigl(
\begin{matrix} 1 &0 \\y_2  & 1
\end{matrix}\Bigr)\,
\Bigl(
\begin{matrix} 1 & -\l x_2 \\  0& 1
\end{matrix}\Bigr)\,
\Bigl(
\begin{matrix}1  & \frac{\l^2x_2}{1+\l }\\ 0  &1
\end{matrix}\Bigr)
= \Bigl(
\begin{matrix} 1+\l & 0\\ 0 & \frac{1}{1+\l }
\end{matrix}\Bigr)\,
\Bigl(
\begin{matrix} 1-\frac{\l E_{00}(2)}{1+\l} & 0\\ 0 &1
\end{matrix}\Bigr).
\end{equation}
This can be checked by direct multiplication using the equalities
$y_2x_2=1$, $x_2y_2= 1-E_{00}(2)$, $y_2E_{00}(2) =0$ and
$E_{00}(2)x_2=0$  in the algebra $\mS_1(2)$. Since the first six
matrices in the equality belong to the group $E_\infty
(\mS_1(2))$, the last matrix $c=\Bigl(
\begin{matrix} 1-\frac{\l E_{00}(2)}{1+\l} & 0\\ 0 &1
\end{matrix}\Bigr)$ belongs to the group $E_\infty (\mS_1(2))$ as
well and can be written as
$$ c=E_{00}(1)(1-\frac{\l }{1+\l}E_{00}(2))+1-E_{00}(1)=
1-\frac{\l }{1+\l}E_{00}(1)E_{00}(2)\in (1+F_2)^*.$$ Since the map
$\v :K\backslash \{ -1\}\ra  K\backslash \{ -1\}$, $ \l \mapsto
-\frac{\l }{1+\l}$, is a bijection  ($\v^{-1} = \v$), all the
elements $b$ belong to the group $E_\infty (\mS_1(2))$. The proof
of the proposition is complete. $\Box$

$\noindent $

{\bf Proof of Theorem \ref {aInt24Apr9}}. The third statement
follows at once from the second one. To prove the remaining two
statements we use induction on $n$. The case when $n=1$ is Theorem
4.6, \cite{shrekaut}.  So, let $n>1$ and we assume that the first
two statements hold for all natural numbers $n'<n$. Clearly,
$\mS_n = \bigotimes_{i=1}^n \mS_1(i)$ where $\mS_1(i):= K\langle
x_i, y_i\rangle \simeq \mS_1$. Consider the following ideals of
the algebra $\mS_n$:
$$\gp_1:=F\t \mS_{n-1},\; \gp_2:= \mS_1\t F\t \mS_{n-2}, \ldots ,
 \gp_n:= \mS_{n-1} \t F, \; \ga_n:= \gp_1+\cdots +\gp_n.$$
 Let $K(x_n)$ be the field of fractions of the polynomial algebra
 $K[x_n]$. It follows from the chain of algebra homomorphisms
\begin{equation}\label{SpKn}
\mS_n\ra \mS_n/\gp_n\simeq \mS_{n-1}\t \mS_1/F\simeq \mS_{n-1}\t
K[x_n, x_n^{-1}]\ra \mS_{n-1}\t K(x_n)\simeq \mS_{n-1}(K(x_n))
\end{equation}
and from the induction on $n$ that $\mS_n^*\subseteq
\mS_1(n)+\ga_n = \sum_{i\geq 1}Ky_n^i+K+\sum_{i\geq
1}Kx_n^i+\ga_n$ (since $\mS_{n-1}(K(x_n))^*= K(x_n)^*$, by
induction). By symmetry of the indices $1, \ldots , n$, we have
the inclusion $\mS_n^* \subseteq \bigcap_{j=1}^n(\sum_{i\geq
1}Ky_j^i+K+\sum_{i\geq 1}Kx_j^i+\ga_n) = K+\ga_n$ and so
$$ K^*(1+\ga_n)^*\subseteq \mS_n^* \subseteq \mS_n^*\cap
(K+\ga_n) = K^* \cdot (\mS_n^* \cap (1+\ga_n))= K^* (1+\ga_n)^* =
K^*\times (1+\ga_n)^*$$ since $\ga_n$ is an ideal of the algebra
$\mS_n$ (and so $\mS_n^* \cap (1+\ga_n)= (1+\ga_n)^*$). This
proves statement 1.

It follows from statement 1, (\ref{SpKn}) and induction that
$Z(\mS_n) \subseteq K^* (1+\gp_n)^*$, hence, by symmetry,
$$ Z(\mS_n^*) \subseteq \bigcap_{i=1}^n K^*(1+\gp_i)^*= K^*
(1+\bigcap_{i=1}^n\gp_i)^* = K^* (1+F_n)^*$$ where $F_n
:=\bigcap_{i=1}^n\gp_i$. Since $(1+F_n)^*\simeq \GL_\infty (K)$
(see Section 2, \cite{shrekalg}) and the centre of the group
$\GL_\infty (K)$ is $\{ 1\}$, statement 2 follows. $\Box$


\section{Normal  subgroups of $\GL_\infty (\mS_1)$ and the centres
of the groups  $\GL_\infty (\mS_1)/\SL $ and $E_\infty (\mS_1)/\SL
$ }\label{TGS1K1}

In this section, several normal subgroups of the group $\GL_\infty
(\mS_1)$ are introduced, see (\ref{normsub}) and Propositions
\ref{a7Apr10}.(4). The most important (and non-obvious) is the
normal subgroup $\SL$ (Proposition \ref{a7Apr10}.(4)). The group
$(1+F_2)^*$ is {\em non-canonically} isomorphic to the group
$\GL_\infty (K)$, hence it inherits the determinant homomorphism
$\det $, see (\ref{SLdet}), and $\SL:= \{ u\in (1+F_2)^* \, | \,
\det (u) =1\}$. We will show that the determinant $\det$ and the
group $\SL$ {\em does not} depend on the isomorphism $(1+F_2)^*
\simeq \GL_\infty (K)$. Moreover, the determinant is invariant
under the conjugation of its argument by
 the elements of the group $\GL_\infty (\mS_1)$, Proposition
 \ref{a7Apr10}.(3). This is the central point of this section. It
 implies that  the group $\SL$ is a {\em normal} subgroup of
 $\GL_\infty (\mS_1)$ and is a key fact in  finding  the
 centres of the groups  $\GL_\infty (\mS_1)/\SL $ and
 $E_\infty (\mS_1)/\SL $ (Theorem \ref{A10Apr10}).

In order to prove Theorem \ref{A10Apr10} we use results and
notations of Section \ref{KTG1S}. In particular, we use the
following group isomorphism:
$$ (1+\gp_1 )^* =(1+M_\infty (\mS_1))^*\simeq \GL_\infty (\mS_1),
\;\; 1+\sum a_{ij} E_{ij}(1)\mapsto  1+\sum a_{ij} E_{ij}, $$
where $a_{ij}\in \mS_1= \mS_1(2)$ and $E_{ij}$ are the matrix
units. It is convenient to identify the groups $(1+\gp_1)^*$ and
$\GL_\infty (\mS_1)$ via this isomorphism, i.e. $E_{ij}(1) =
E_{ij}$. Then, by Proposition \ref{c24Apr9},
$$\GL_\infty (\mS_1) = U\ltimes E_\infty (\mS_1)$$
where $U:= \{ \mu (\l ) :=\l E_{00}+1-E_{00}\, | \, \l \in
K^*\}\simeq K^*$, $\mu (\l ) \lra \l$, and the groups $E_\infty
(\mS_1)$ and  $(1+F_2)^*$ are normal subgroups of $\GL_\infty
(\mS_1)$.

As we have seen in the proof of Proposition \ref{10Apr10} the
group $(1+F_2)^*$ is isomorphic to the group $\GL_\infty (K)$.
This isomorphism depends on the choice of the bijection $b$. For
the group $\GL_\infty (K)$, we have the determinant (group
epimorphism) $\det : \GL_\infty (K) \ra K^*$, the short exact
sequence of groups $1\ra \SL_\infty (K) \ra \GL_\infty
(K)\stackrel{\det }{\ra}K^*\ra 1$, and the decomposition
$\GL_\infty (K) = U(K) \ltimes \SL_\infty (K)$ where $U(K)= \{ \mu
(\l ) \, | \, \l \in K^*\}$. Therefore, for the group $(1+F_2)^*$
we have the determinant (group epimorphism) $\det : (1+F_2)^* \ra
K^*$, the short exact sequence of groups 
\begin{equation}\label{SLdet}
1\ra \SL \ra (1+F_2)^* \stackrel{\det }{\ra}K^*\ra 1
\end{equation}
and the decomposition $(1+F_2)^* = U' \ltimes \SL$ where
\begin{eqnarray*}
 U'&:=& \{ \mu' (\l ):=\l E_{00}(1)E_{00}(2)+1-E_{00}(1)E_{00}(2) \, | \,
 \l \in K^*\}\simeq K^*, \; \mu'(\l ) \lra \l ,  \\
 \SL &:=& \{ u\in (1+F_2)^* \, | \, \det (u)=1\}.
\end{eqnarray*}
The group $\SL$ is generated by the elements $1+\l
E_{\alpha\beta}$ where $\l \in K$ and $\alpha, \beta \in \N^2$
such that $\alpha \neq \beta$.  Theorem 8.1, \cite{shrekaut}, says
that the map $\det : (1+F_2)^* \ra K^*$ {\em does not} depend on
the choice of the bijection $b$.

\begin{theorem}\label{6Mar9}
{\rm (Theorem 8.1, \cite{shrekaut})}  Let $\CV = \{ V_i\}_{i\in
\N}$ be a finite dimensional vector space filtration on $P_2$
(i.e. $V_0\subseteq V_1\subseteq \cdots$ and $P_2= \cup_{i\in \N}
V_i$) and $a\in (1+F_2)^*$. Then $a(V_i) \subseteq V_i$ and $\det
(a|_{V_i}) = \det (a|_{V_j})$ for all $i,j\gg 0$. Moreover, this
common value of the determinants does not depend on the filtration
$\CV$ and, therefore, coincides with the determinant in
(\ref{SLdet}).
\end{theorem}
By Theorem \ref{6Mar9},  the group $\SL$ does not depend on the
choice of the bijection $b$. The algebra $\mS_2$ admits the {\em
involution}:
$$ \eta : \mS_2\ra \mS_2, \;\; x_i\mapsto y_i, \;\; y_i\mapsto
x_i, \;\; i=1,2,$$ i.e. it is a $K$-algebra {\em anti-isomorphism}
 ($\eta (ab) = \eta (b) \eta (a)$ for all elements $a,b\in
 \mS_2$) such that $\eta^2= {\rm id}_{\mS_2}$, the identity map on
 $\mS_2$. It follows that
\begin{equation}\label{etaEij}
 \eta (E_{ij} (k)) = E_{ji}(k))\;\; {\rm and}\;\; \eta
 (E_{\alpha\beta}) = E_{\beta \alpha}
\end{equation}
for all elements $i,j\in \N$, $k=1,2$, and $\alpha , \beta \in
 \N^2$. Therefore, $\eta (\gp_k) = \gp_k$ and $\eta (F_2) = F_2$.
 It is easy to see that $\eta ((1+\gp_k)^*)= (1+\gp_k)^*$,
 $\eta ((1+F_2)^*)= (1+F_2)^*$ and $\eta (E_\infty (\mS_1(k)))=
 E_\infty (\mS_1(k))$.

 The polynomial algebra $P_2$ is equipped with the {\em cubic}
filtration $\CC := \{ \CC_m:=\sum_{\alpha\in C_m}Kx^\alpha\}_{m\in
\N}$ where $C_m:= \{ \alpha \in \N^2\, | \, {\rm all} \;
\alpha_i\leq m\}$. The filtration $\CC$ is an ascending, finite
dimensional filtration such that $P_2= \bigcup_{m\geq 0}\CC_m$ and
$\CC_m\CC_l\subseteq \CC_{m+l}$ for all $m,l\geq 0$.

\begin{proposition}\label{a7Apr10}
\begin{enumerate}
\item For all elements $a\in (1+F_2)^*$, $\det (\eta (a))=\det
(a)$. \item Let $ a\in \GL_\infty (\mS_1)=(1+\gp_1)^*$ and let
$\CV = \{ V_i\}_{i\in \N}$ be an ascending finite dimensional
filtration on $P_2$ such that $a(V_i) \subseteq V_i$ for all $i\gg
0$. Then $\det (aba^{-1}) = \det (a)$ for all elements $b\in
(1+F_2)^*$. \item For all elements $a\in \GL_\infty (\mS_1)$ and
$b\in (1+F_2)^*$, $\det (aba^{-1}) = \det (a)$. \item The group
$\SL$ is a normal subgroup of $\GL_\infty (\mS_1)$. Moreover, each
subgroup $N$ of the group $(1+F_2)^*$ that contains the group
$\SL$ is a normal subgroup of $\GL_\infty (\mS_1)$.
\end{enumerate}
\end{proposition}

{\it Proof}. 1. Recall that $(1+F_2)^* = U'\ltimes \SL$. It is
obvious that $\eta (u) = u$ for all elements $u\in U'$ and $\eta
(\SL) \subseteq \SL$ (by (\ref{etaEij}), the set of standard
generators of the group $\SL$ is invariant under the action of the
involution $\eta$). An element $a\in (1+F_2)^*$ is a unique
product $a=us$ for some elements $u\in U'$ and $s\in \SL$. Now,
$$ \det (\eta (a)) = \det (\eta (us))= \det (\eta (s) u) = \det
(u) = \det (a).$$

2. Applying Theorem \ref{6Mar9} to the element $aba^{-1} \in
(1+F_2)^*$ we have the result, for all  $i\gg 0$:
$$ \det (aba^{-1}) = \det (aba^{-1}|_{V_i})=
 \det (a|_{V_i}\cdot  b|_{V_i} \cdot a^{-1}|_{V_i})= \det (b|_{V_i}) = \det
 (b).$$

 3. We deduce statement 3 from the previous two. Recall that we
 have
 identified the groups $(1+\gp_1)^*$ and $\GL_\infty (\mS_1)$. By
 Theorem \ref{d24Apr9}, $\GL (\mS_1) = U\ltimes E_\infty
 (\mS_1)$. Since the subgroup $(1+F_2)^*$ of $\GL_\infty (\mS_1)$ is normal,
 it suffices to show  that statement 3 holds for generators of
 the groups $U$ and $E_\infty (\mS_1)$. Since $U=\{ \mu (\l ) \, |
 \, \l \in K^*\}$ and $\mu (\l ) (\CC_m) \subseteq \CC_m$ for all
 $m\geq 0$ where $\CC =\{ \CC_m\}_{m\in \N}$ is the cubic
 filtration on the polynomial algebra $P_2$, we see that $\det
 (\mu (\l ) b\mu (\l )^{-1}) = \det (b)$, by statement 2. The
 group $E_\infty (\mS_1)$ is generated by the elements (where $\l
 \in K$, $i\neq j$): $a_{ij}= 1+E_{ij}\l x_2^n$ where $n\geq 1$;
 $b_{ij}= 1+E_{ij}\l y_s^m$ where $m\geq 0$; and $1+E_{ij}f$ where
 $f\in F_2$. Statement 3 holds for the elements $1+E_{ij}f$ since
 they belong to the group $(1+F_2)^*$.

 Since $b_{ij}(\CC_k)\subseteq \CC_k$ for all $k\gg 0$, statement
 3 holds for the elements $b_{ij}$, by statement 2. Since $\eta
 (a_{ij}) = b_{ji}$, statement 3 holds for the elements $a_{ij}$,
 by statements 1 and 2. In more detail,
 $$ \det (a_{ij}ba_{ij}^{-1}) = \det (\eta (a_{ij}ba_{ij}^{-1}))
 =\det (b_{ji}^{-1} \eta (b) b_{ji}) = \det (\eta (b)) = \det
 (b). $$
4. By statement 3, the group $\SL$ is a normal subgroup of
$\GL_\infty (\mS_1)$. If $N$ is  a subgroup of $(1+F_2)^*$
containing $\SL$ then $N=N'\ltimes \SL$ for a subgroup $N'=\{ \mu
(\l ) \, | \, \l \in K'\}$ of $U'$ where $K'$ is an additive
subgroup of the field $K$. Since $K'=\det (N)$ and $\det
(aNa^{-1}) = \det (N)=K'$, by statement 3, we see that $aNa^{-1}
\subseteq N$, i.e. $N$ is a normal subgroup of $\GL_\infty
(\mS_1)$. $\Box $

By Proposition \ref{a7Apr10}.(4), there is a chain of normal
subgroups of the group $\GL_\infty (\mS_1)$: 
\begin{equation}\label{normsub}
\SL\subset (1+F_2)^* \subset E_\infty (\mS_1) \subset \GL_\infty
(\mS_1).
\end{equation}
Using the fact that $(1+F_2)^* = U'\ltimes \SL$, we have the chain
of normal subgroups of the factor group $\GL_\infty (\mS_1) /
\SL$: $$U'\subset E_\infty (\mS_1)/\SL \subset \GL_\infty (\mS_1)
/ \SL.$$
\begin{theorem}\label{A10Apr10}
The group $U'\simeq K^*$ is the centre of both groups $\GL_\infty
(\mS_1) / \SL$ and $E_\infty (\mS_1)/\SL $.
\end{theorem}

{\it Proof}. By Proposition \ref{a7Apr10}.(3), the group $U'$
belongs to the centres of both groups. The algebra homomorphism
$\mS_1\ra \mS_1/F\simeq K[x_2, x_2^{-1}]$ yields the group
homomorphisms:  $\GL_\infty (\mS_1)\stackrel{\v}{\ra} \GL_\infty
(K[x_2, x_2^{-1}])$ and $E_\infty (\mS_1)\stackrel{\psi}{\ra}
E_\infty (K[x_2, x_2^{-1}])$. By Proposition \ref{c24Apr9}, $\im
(\v ) = U(K)\ltimes E_\infty (K[x_2, x_2^{-1}])$ and $\im (\psi )
= E_\infty (K[x_2, x_2^{-1}])$. Since the groups $ \im (\v ) $ and
$\im (\psi )$ have trivial centre and $\ker (\v ) = \ker (\psi ) =
(1+F_2)^* / \SL = U'\ltimes \SL / \SL \simeq U'$,  the result
follows. $\Box $

$\noindent $

$\noindent $

Department of Pure Mathematics

University of Sheffield

Hicks Building

Sheffield S3 7RH

UK

email: v.bavula@sheffield.ac.uk

\end{document}